\documentclass[11pt]{article}

\usepackage{amsfonts,amsmath,amssymb,amsthm,enumerate,graphicx,pgfplots,hyperref,float,url}
\def\ballk{\mathcal B_{\theta}}
\def\bsup{\sup\left\{ }
\def\binf{\inf\left\{ }
\def\bset{\{ }
\def\call{\textrm{call}}

\def\cpl{\textrm{Cpl}}

\def\dc{d_c}

\def\FF{\mathcal{F}}

\def\half{\frac{1}{2}}
\def\intrd{\int_{\RR^d}}

\def\dnu{\,d\nu}
\def\dnuu{\,d\nu_0}
\def\mina{\min_{\vmin\leq \alpha\leq \vmax}}	
\def\infl{\inf_{\lambda\geq 0}}
\def\nun{\nu_0}
\def\NN{\mathcal N}
\def\NNN{\mathbb{N}}
\def\RR{\mathbb{R}}
\def\intqb{\int_{\qb}^{\infty}}

\def\rrbeta{\frac{1}{1-\beta}}

\def\einf{\right\}}
\def\eset{\}}

\def\st{:\ }
\def\suchthat{\st}
\def\esup{\right\}}
\def\cvar{\textrm{ES}_{\beta}}
\def\varb{\textrm{VaR}_{\beta}}
\def\varnb{\textrm{VaR}_{\beta}^{\nu}}

\def\cvarnb{\cvar^{\ \nu}}
\def\cvarnnb{\cvar^{\ \nu_0}}

\def\flc{f^{\lambda c}}

\def\mi{m_i}
\def\mix{\langle m_i,x\rangle}
\def\nn{\mathcal N}

\def\nuu{\nu_0}
\def\PP{\mathcal P}
\def\qb{q_{\beta}}
\def\rlambda{\frac{1}{\lambda}}
\def\rtlambda{\frac{1}{2\lambda}}
\def\supin{\sup_{i\in\NNN}}
\def\supyrd{\sup_{y\in\RR^d}}
\def\tr{\tilde\rho}

\def\vmin{v_{min}}
\def\vmax{v_{max}}

\def\nin{\nu\in\nn}

\def\var{VaR}

\newtheorem{theorem}{Theorem}

\theoremstyle{definition}
\newtheorem{example}[theorem]{Example}
\newtheorem{definition}[theorem]{Definition}

\title{Distributionally robust Expected Shortfall for convex payoffs}

\date{}

\author{ Gusti van Zyl}

\begin{document}

\maketitle

\noindent{\bf Correspondence:} Gusti van Zyl, 
Department of Mathematics and Applied Mathematics, University of Pretoria,
South Africa.

\noindent Email: gusti.vanzyl@up.ac.za
\medskip

\noindent{\bf Funding:} Research supported in part by the National Research Foundation of South
Africa (Grant Number 146018).
\medskip

\noindent{\bf Conflict of interest:} The author has no conflict of interest to declare.\medskip

%\noindent {\bf Data availability:} Data sharing not applicable to this article as no datasets were generated or analysed during the current study.

%\noindent {\bf Keywords:} {\it optimal transport; $\lambda$c-transform; quadratic cost function; robust risk measure; coherent risk measure.}

\begin{abstract}
We study distributionally robust expected shortfall when the distribution of the underlying is perturbed by a size quantified with optimal transport distance based on the quadratic cost function. In the dual version of the robust expectation problem, which is part of the robest expected shortfall problem, the computation of the so-called $\lambda c$-transform $\flc$ of payoff $f$ is required. We show that under the quadratic cost function there exists a tractable representation of $\flc$, if $f$ is convex. Furthermore, we show that robust expected shortfall can be characterized as the solution of a 2-dimensional minimization problem. We apply these results to obtain a closed-form formula for robust, with respect to the risk-neutral distribution, expected shortfall of an unhedged call option, from the point of view of the writer. %, as well as that of a portfolio mixing underlying shares with a call and a put option. 
\end{abstract}

\section{Introduction}
The study of distributional risk attempts to quantify the ``consequences of using the wrong models [or statistical distributions] in risk measurement, pricing and portfolio selection" \cite{BreuerCsiszár2016}. Causes of distributional risk include bad choice of stochastic process, wrong dependence assumptions, and parameter drift. Examples of wrong distributional assumptions are well known. The demise of LCTM in 1999 is partly blamed on overreliance on Gaussian logreturns, and the breakdown of credit models in the 2008 financial crisis, on the unexpected behaviour of correlations during crisis periods. Model risk \cite{morini2011understanding}  is broader than distributional. It may also be mentioned that the study of distributional risk is intended here not as a substitute for the correct dynamics, but as an attempt to probe, and possibly compare, the distributional risks of different portfolios, if no further knowledge about the distribution is assumed than that it may be off from the correct distribution by a certain distributional ``distance." 

One could assess distributional risk using parametric or non-parametric approaches. Parametric approaches include describing the sensitivity of an option price with regard to volatility and other statistically important parameters. A related approach is to parameter intervals.

Non-parametric approaches are usually based on the question of how much a financially relevant calculation could be affected if the underlying ``baseline" or nominal distribution is varied within an ``uncertainty set" of alternative distributions, the set of alternatives not characterized by a few parameters. Typically the uncertainty set comprises those distributions that differ from the baseline by at most a certain amount $\theta>0$. There are several prominent ways to quantify this difference between the baseline and alternative distributions. In quantitative finance, Glasserman and Xu \cite{glasserman_robust_2014} pioneered the use of relative entropy as distance measure between distributions. See, for example Feng et al \cite{feng2018quantifyingmodelriskinherent} for an application to option pricing models. Another quantification, which is the focus of this paper, is by  optimal transport-like distances of measures. Other statistical ``divergence" measures, for example those suited to the tails of distributions, are also possible.

The above-mentioned list of categories is of course by no means complete. It may also be mentioned that the term ``distributionally robust", which is widespread in the literature, is not meant to suggest robustness with respect to all possible changes in distribution, which is in general an impossible goal. Instead, in the non-parametric approaches referred to above, the distributional robustness is always relative to an ambiguity tolerance parameter $\theta$. For larger values of $\theta$ the calculation is more conservative. The non-robust case corresponds to $\theta=0$.

The optimal transport distance $d_c$ to quantify the distance on the set $P(\RR^d)$ of probability measures on $\RR^d$, is based on a cost function $c(x,y)$ for $x,y\in\RR^d$. In this paper $c(x,y)=\frac{1}{2}\|x-y\|^2$, called the quadratic cost function, with $\|\|$ denoting the Euclidean norm. The robust expectation problem is 
\begin{equation}\label{eq:primal}
P:=\sup\left\{\int fd\nu:\ \nu \in \NN\right\},	
\end{equation}
where $\nn=\{\nu\in P(\RR^d):\ \dc(\nu,\nu_0)\leq \theta \}$.  In this optimization problem the uncertainty set, involving a set of probability distributions or measures, is infinite dimensional. For applications a duality result such as that by Bartl et al \cite[Theorem 2.4 and Equation (7)]{bartl_computational_2019}
\begin{equation}\label{eq:dual-robust}
	P=D:=\inf_{\lambda\geq 0}\left\{\lambda\theta+\int \flc d\nu_0\right\},
\end{equation}
where $\flc(x)=\supyrd \{f(y)-\lambda c(x,y)\}$, is usually needed. 
Thus the infinite-dimensional problem is reduced to one-dimensional optimization problem, provided that $\flc$ is known. The main observation of this paper, which we show in section \ref{se:piecewise}, is that there exists a tractable representation of $\flc$ if $f$ is convex. This case is important in mathematical finance since non-negative  combinations of call and put options are convex. 

In section \ref{se:cvar-call} we use the dual problem for robust expectation to derive a representation of robust Expected Shortfall in terms of an optimization over two parameters, and apply this to derive an analytical expression for the robust Expected Shortfall of a call option, from the point of view of the writer, under the risk-neutral distribution.
%In section \ref{se:portfolio} we apply this representation to portfolio robust expected shortfall minimization. 

\section{$\lambda c$-transform for convex, piecewise-linear, functions}\label{se:piecewise}
Consider a probability space $(\Omega,\FF,P)$. For any random variable $Y:\Omega\to\RR^d$, its distribution $\mu_Y$ defined by $\mu_Y(A)=P(Y\in A)$ is a real-valued probability measure on the Borel subsets of $\RR^d$. Therefore we will be interested in probability measures, referred to as distributions, on $X=\RR^d$. 

We recall the optimal transport distance \cite{villani_optimal_2009} between Borel probability measures $\mu$ and $\nu$, 
\begin{equation}\label{eq:dc}
	d_c(\mu,\nu):=\binf \int_{X\times X} c(x,y)d\pi(x,y) \st \pi\in \cpl(\mu,\nu)\einf,
\end{equation} where $\cpl(\mu,\nu)$ denotes the set of ``couplings" of $\mu$ and $\nu$; that is, the set of probability measures $\pi$ on $X\times X$ with first marginal distribution equal to $\mu$, and second to $\nu$. A coupling will always exist, take for example $\pi$ as the product measure $\mu\otimes\nu$.   With the quadratic cost function, it holds \cite[Theorem 4.1]{villani_optimal_2009} that $d_c(\mu,\nu)<\infty$ if $\mu$ and $\nu$ have finite second moments; that is, if $\int_X \|x_0-x\|^2d\mu(x)<\infty$ and $\int_X \|x_0-x\|^2d\nu(x)<\infty$ for any $x_0\in X$. 

In the case of the quadratic cost function it is well known the infimum in Equation (\ref{eq:dc}) is obtained by a transport map $T:\RR^d\to\RR^d$ that defines $\nu=T_\#(\mu)$, in other words the image measure induced by $\nu_0$ acting through $T$.

 To give some intuition, we consider a special case of discrete distributions $\mu$ and $\nu$, $X=\RR$, see Figure \ref{fig:Wasserstein}. 
\begin{figure}[H]
\begin{tikzpicture}[scale=1]
	% Axes
	\draw[->] (-1,0) -- (3.2,0) node[right] {$x$};
	\draw[->] (0,0) -- (0,1.5) node[above] {$p(x)$};
	% Distribution 1.
	\draw[thick, blue] (0,0) -- (0,0.25);
	\filldraw[blue] (0,0.25) circle (2pt);
	\node[above left] at (0,0.25) {$\frac 14$};
	\draw[thick, blue] (0.5,0) -- (0.5,0.75);
	\filldraw[blue] (0.5,0.75) circle (2pt);
	\node[above left] at (0.5,0.75) {$\frac 34$};
	% Label
	\node[below] at (0,-0.1) {0};
	\node[below] at (0.5,-0.1) {$\frac 12$};
	% Distribution 2.
	\draw[thick, red] (2,0) -- (2,0.5);
	\filldraw[red] (2,0.5) circle (2pt);
	\node[above left] at (2,0.5) {$\frac12$};
	\draw[thick, red] (3,0) -- (3,0.5);
	\filldraw[red] (3,0.5) circle (2pt);
	\node[above left] at (3,0.5) {$\frac12$};
	% Label
	\node[below] at (0,-0.1) {0};
	\node[below] at (0.5,-0.1) {$\frac 12$};
	\node[below] at (2,-0.1) {$2$};
	\node[below] at (3,-0.1) {$3$};
	% Legends (optional)
	\node[blue] at (1,-1) {$\mu(\{0\})=\frac 14,\ \mu(\{\frac 12\})=\frac 34$,};
	\node[red] at (1,-1.5) {$\nu(2)=\nu(3)=\frac 12$.};
%	\node at (6,1) {$W_2(\mu,\nu)=\frac12 |3-0.5|^2+\frac14 |2-0.5|^2$};
%	\node at (6,0.5) {$+\frac14|2-0|^2=4.6875$};
\end{tikzpicture}	
\caption{Distributions to illustrate optimal transport distance under quadratic cost}	
\label{fig:Wasserstein}
\end{figure}

One way to transport $\mu$ to $\nu$ is to move $\frac 12$ probability mass from $x=\frac 12$ to $x=3$, and all the remaining masses to $x=2$. The average cost of this transport is $\frac 12 c(\frac 12,3)+\frac 14 c(\frac 12,2)+\frac 14 c(0,2)=2.34375.$ This turns out to be the most efficient ``transport plan" for these distributions, more efficient for example than moving the mass at $x=0$ to $x=3$. Therefore $d_c(\mu,\nu)=2.34375$. Optimal transport distance for distributions that are not necessarily discrete, can be defined using an approximation procedure.

As mentioned in the introduction, the robust expected value problem is tractable once the $\lambda c$-transform $\flc$ is available. In the case of convex $f$, the following representation of $\flc$ is known.
\begin{theorem}\label{thm:lc}
	Let $f:\RR^d\to\RR$ be a convex function, represented as $f(x)=\supin\{\mix+c_i\},$ with $\mi\in\RR^d,c_i\in\RR,i\in\NNN.$ Let $c$ denote the cost function $c(x,y)=\half \|x-y\|^2$.
	Then $$\flc(x)=\supin \left\{\mix+c_i+\frac{\|m_i\|^2}{2\lambda}\right\}.$$
	In particular, if $f$ is piecewise-linear and convex, then so is $\flc$.
\end{theorem}
\begin{proof} It is well known that every convex function on $\RR^d$ is a pointwise supremum of a countable collection of affine functions. 
Let $f_i(x)=\mix+c_i$. Then 
\begin{eqnarray*}
	\flc(x) &=& \supyrd \{f(y)-\lambda c(x,y)\} \\
	&=& \supyrd \{\supin f(y)-\lambda c(x,y)\} \\
	&=& \supin \flc_i(x).	
\end{eqnarray*}
To find $\flc_i$ observe that the mapping $y\mapsto f_i(y)-\rtlambda \|x-y\|^2$ is concave and differentiable, so it is sufficient to apply the first-order condition $m_i-\lambda(x-y)=0$; that is, $y=x+\frac{1}{\lambda}m_i$. 
So $\flc_i(x)=f_i(x+\rlambda m_i)=\mix +c_i+\rtlambda\|\mi\|^2$. 
\end{proof}
This result shows that under quadratic transportation cost, the $\lambda c$-transform shifts each affine component upwards by a constant proportional to the square of its slope. A partial heuristic interpretation is given in Figure \ref{fig:lc}, where the $\lambda c$-transforms, of payoff functions that differ by a constant $k$, also differ by $k$. This is in line with the intuition that a payoff, independent of the underlying, is not exposed to the distribution risk of the underlying.

\begin{figure}[H]
\begin{tikzpicture}
	% Axes
	\draw[->] (-1.5, 0) -- (6, 0) node[right] {$x$};
	\draw[->] (0, 0) -- (0, 3.0) node[above] {Payoff};
	\draw[-] (3,0.1) -- (3, -0.1) node[below] {$K$};
	% Payoff function 1
	\draw[-, color=blue] (-1.5, 0) -- (3, 0);
	\draw[->, color=blue] (3, 0) -- (4, 1) node[right] {$f(x)$};
	% Payoff function 2
	\draw[-, color=green] (-1.5, 0) -- (3, 0) ;
	\draw[->, color=green] (3, 0) -- (4, 2) node[right]{$f_1(x)$} ;
	% Payoff function 3
	\draw[-, color=purple] (-1.5, 2) -- (3, 2);
	\draw[->, color=purple] (3, 2) -- (4, 3) node[above]{$f_2(x)$\qquad};
	% Normal curve centered at x=1
	\draw[domain=-1:5, smooth, samples=100, thick]
	plot(\x, {3*0.3989 * exp(-((\x - 1)^2)/2)});
	\draw[domain=-1:5, smooth, samples=100, thick, color=red]
	plot(\x, {3*0.325 * exp(-((\x - 1)^2)/4.5)});
\end{tikzpicture}
\caption{A heuristic interpretation of Theorem \ref{thm:lc}. Slope matters but vertical intercept does not. In comparison with $f$, payoff $f_1$ is  more sensitive to the change of distribution, but $f_2$ not.}
\label{fig:lc}
\end{figure}	

\begin{example}
	Consider a call option payoff $f(x)=\max\{x-K,0\}$. Then $$\flc(x)=\max\{x-K+\frac{1}{2\lambda},0\}=\max\{x-(K-\frac{1}{2\lambda})\},$$ in other words the effect of this part of the robustification process is to reduce the strike price by $\frac{1}{2\lambda}$. 
	
	This can be considered a special case of \cite[Example 2.14]{bartl_computational_2019} which was calculated by direct use of first-order conditions, instead of Theorem \ref{thm:lc}. Consider $$h^{\lambda c,\alpha}(x):=\displaystyle\sup_{y\in \RR}\left((y-k)^++\alpha(y-s)-\frac{\lambda}{2} (y-x)^2\right)$$ is to be computed for $\lambda>0$. In our notation this is $\flc(x)$ where $$f(x)=(x-k)^++\alpha(x-s)=\max\{\alpha y-\alpha s,(\alpha+1)y-k-\alpha s\}.$$ By Theorem \ref{thm:lc}, 
\begin{eqnarray*}
	\flc(x)&=& \max\left\{\alpha x-\alpha s+\frac{\alpha^2}{2\lambda},(\alpha+1)x-k-\alpha s+\frac{(\alpha+1)^2}{2\lambda}\right\},
\end{eqnarray*}
from which it is easy to recover the result \cite[Example 2.14]{bartl_computational_2019}
$$h^{\lambda c,\alpha}(x)=\left(x-\left(k-\frac{2\alpha+1}{2\lambda}\right)\right)^++\alpha(x-s)+\frac{\alpha^2}{2\lambda}.$$
\end{example}

\begin{example} Put option with strike $K_1$ and two call options with strike $K_2>K1$. 
	Then $f(x)=\max\{K_1-x,0\}+2\max\{x-K,0\}=\max\{-x+K_1,2x-2K_2,0\}$.
Using Theorem \ref{thm:lc}, 
$$\flc(x)=\max\left\{-x+K_1+\frac{1}{2\lambda},2x-2K_2+\frac{2}{\lambda},0\right\}=\max\left\{-x+K_1+\rtlambda,2(x-(K_2-\frac{1}{\lambda})),0\right\}.$$	\medskip

\begin{figure}[H]
\begin{tikzpicture}
	%	% Axes
	\draw[->] (-1, 0) -- (3, 0) node[right] {$x$};
	\draw[->] (0, 0) -- (0, 3.0) node[above] {Payoff};
	%\draw[-] (3,0.1) -- (2, -0.1) node[below] {$K$};
	%	
	%	% Payoff function 1 
	\def\strike{1}
	\def\strikePut{1} \def\strikeCall{2} \def\rectwolambda{0.25} 
	\def\frtlambda{0.5} % 2x0.25
	\draw[-, color=blue] (0, \strikePut) -- (\strikePut,0);
	\draw[->, color=blue] (\strikeCall, 0) -- (\strikeCall+1, 2) node[right] {$f(x)$};
	%	% Payoff function 2
	%\draw[-, color=brown,dashed] (0, \strike+\rectwolambda) -- (\strike+\rectwolambda,0);
	%\draw[-, color=brown, dashed] (\strike-\rectwolambda,0) -- (\strike-\rectwolambda+2,2);
	\draw[-,color=brown] (0,\strikePut+\rectwolambda) -- (\strikePut+\rectwolambda,0);
	\draw[->,color=brown] (\strikeCall-\frtlambda,0) -- (\strikeCall-\frtlambda+1,2) node[left] {$\flc(x)$}; ;%(\strikeCall+4\rectwolambda+1, 4\rectwolambda+2) 
	% Legend
	\node at (\strike,-0.3) {$K$};
\end{tikzpicture}
\caption{$\lambda c$-transform adjusts $f$ based on the slope of the affine component.}
\end{figure}
\end{example}

\begin{example}\label{ex:three-asset}
We consider a skewed three-asset portfolio similar to the one used by \cite{bertsimas_shortfall_2004} to 
	 demonstrate portfolio optimization against ``shortfall" as risk metric. 
%Shares independent LogNormal($\mu_i,\sigma_i$)  chosen so A,B,C have same means and variances; $r_f=2.5\%$.
Asset A comprises a share, asset B an independent underlying share and 0.75 at-the-money (ATM) call options on the underlying, and asset C an independent underlying share and 0.75 ATM put options on the underlying. The risk-free rate is $2.5\%$. The payoffs are
	\begin{eqnarray*}
		f_A(x_1)&=& x_1,\\ 
		f_B(x_2)&=& x_2+0.75\max\{x_2-K_2,0\}-0.75(1.025)C_0, \\
		f_C(x_3)&=& x_3+0.75\max\{K_3-x_3,0\}-0.75(1.025)P_0,	
	\end{eqnarray*}
	where $C_0$ (resp. $P_0$) is the price of a call (resp. put) option on the underlying.
	The map $$x=\begin{bmatrix}
		x_1 \\ x_2 \\ x_3
	\end{bmatrix}\mapsto f(x):=w_1f_A(x_1)+w_2 f_B(x_2)+w_3f_C(x_3)$$ is a sum of convex functions, hence, convex, if weights satisfy $w_1,w_2,w_3\geq 0$.

	For the use of Theorem \ref{thm:lc} we rewrite this as 
	\begin{equation}\label{eq:three-asset-payoff}
		f(x)=\max_{i=1,\dots,4} \left\{\langle m_i,x\rangle +c_i\right\},
	\end{equation}
	where 
	\begin{eqnarray*}
		m_1&=& \begin{bmatrix}
			w_1 \\ 1.75 w_2 \\ 0.25 w_3
		\end{bmatrix}, 
		m_2= \begin{bmatrix}
			w_1 \\ 1.75 w_2 \\ w_3
		\end{bmatrix},
		m_3=\begin{bmatrix}
			w_1 \\ w_2 \\ 0.25 w_3
		\end{bmatrix},
		m_4=\begin{bmatrix}
			w_1 \\ w_2 \\ w_3
		\end{bmatrix}
	\end{eqnarray*}	
	and \begin{eqnarray*}
	c_1&=&0.75 w_2 (-k- C_0 (1+r))+0.75 w_3 (k-P_0 (1+r)),\\
	c_2&=&0.75  w_2 (-k- C_0 (1+r))+0.75  w_3 (-P_0 (1+r)),\\
	c_3&=&0.75  w_2 (- C_0 (1+r))+0.75  w_3 (k-P_0 (1+r)),\\
	c_4&=&0.75  w_2 (- C_0 (1+r))+0.75  w_3 (-P_0 (1+r)).
	\end{eqnarray*}
	(The four terms correspond to the four possibilities on whether the put- and/or call option is in- or out of the money.)
	Using Equation (\ref{eq:three-asset-payoff}) and Theorem \ref{thm:lc} we have
	\begin{equation}\label{eq:flc-three-asset} 
		\flc(x)=\max_{i=1,\dots,4} \left\{\langle m_i,x\rangle +c_i+\frac{\|m_i\|^2}{2\lambda}\right\}.
	\end{equation}
\end{example}

\section{Robust Expected Shortfall}\label{se:cvar-call}
We consider expected shortfall as a risk measure $\rho$ on payoff functions on $\RR^d$. For risk measures on functions, see \cite{delbaen_monetary_2024}. In the distributional robustification process, we associate a risk measure $\rho_{\nu}$ to each $\nu$ in the uncertainty set. In particular, Expected Shortfall of the same payoff under different probability distributions, are seen as different risk measures in the sense of giving different risk assessments. 
\begin{definition}
Let $\NN$ be a set of Borel probability measures on $\RR^d$ such that $\nu_0\in A$, and $\rho_{\nu}$ a risk measure based on $\nu$, for every $\nu\in \NN$. Then $\tilde\rho(f):=\sup\{\rho_{\nu}(f):\ \nu\in \NN\}$ is the robustified risk measure. 
\end{definition}
It may be observed that if each $\rho_{\nu}$ is a coherent risk measure, then so is $\rho$. 

For simplicity we consider a liability $f$ on an atomless space, so that $\cvarnb(f)=\rrbeta \int_{w}^{\infty} f(x)d\nu(x)$, where $w:=\varnb(f)$ is the $\beta$-quantile of $f$ under the measure $\nu$, and $\beta$ the confidence level under which Expected Shortfall is computed. Let us denote the baseline risk-neutral measure by $\nu_0$. As in \cite{bartl_computational_2019} we use the ``minimizing property" that characterizes $\cvar$, see Rockafellar and Uryasev \cite[Theorem 1]{rockafellar_optimization_2000}, 
\begin{equation}\label{eq:rockafellar}
	\cvarnb(f)=\min_{\alpha\in\RR} \left\{\alpha+\rrbeta\int_{\RR^d} (f(x)-\alpha)^+\dnu(x)\right\}.
\end{equation}
This expression will be referred to as the dual formula for Expected Shortfall, as it can also be derived using the dual representation of that risk measure.
We recall from \cite{rockafellar_optimization_2000} that $\alpha^*:=\varb(f)$ is a minimizer in this expression. In this section $d=1$. Our uncertainty set is
\begin{equation*}
	\nn:=\ballk(\nun):=\bset \nu\in \PP \suchthat d_c(\nu,\nun)\leq \theta\eset.
\end{equation*}

\begin{theorem} \label{th:robust-ES}
	Let $f:\RR^d \to\RR$ be Lipschitz continuous. The robust $\cvar(f)$ namely 
	\begin{equation*}
	\tr(f):=\displaystyle\sup_{\nin}\rho_{\nu}(f)=\bsup \cvarnb(f) \st d_c(\nu,\nu_0)\leq \theta \esup,	
	\end{equation*}
	 where $c$ denotes the quadratic cost function, 
	satisfies 
	\begin{equation}\label{eq:inf-alpha-lambda}
		\tr(f)=\binf \alpha+\frac{\lambda\theta}{1-\beta}+	
		\rrbeta \intrd ((f(x)-\alpha)^+)^{\lambda c}\ d\nun(x)
		\suchthat \alpha\in \RR, \lambda\geq 0 \einf.
	\end{equation}
\end{theorem}

\begin{proof} Using the definition of $\tilde\rho$ and Equation (\ref{eq:rockafellar}),
	\begin{eqnarray*}
		\tr(f) &=& \sup_{\nin} \cvarnb(f)\\
		&=& \sup_{\nin}\min_{\alpha \in\RR} \left\{\alpha+\rrbeta\intrd (f(x)-\alpha)^+\dnu(x)\right\}\\
		&=& \sup_{\nin}\min_{\vmin\leq \alpha\leq \vmax} \left\{\alpha+\rrbeta\intrd (f(x)-\alpha)^+\dnu(x)\right\},
	\end{eqnarray*}	where $\vmin$ and $\vmax$ are lower and upper bounds, respectively, of $\varnb(f)$ over $\nin$. 
	
	\label{upper-bound-proof}
		To see that an upper bound $\vmax$ exists, let $\|f\|_{Lip}$ be the Lipschitz constant of $f$ and $W_1,W_2=\sqrt{d_c}$ denote the respective p-Wasserstein distances. It follows from the Markov inequality, 1-Wasserstein duality \cite[Remark 6.5]{villani_optimal_2009}, and the basic property $W_1\leq W_2$ \cite[Remark 6.6]{villani_optimal_2009} that for $t>0$ ,
		\begin{eqnarray*}
		\nu(\{x:f(x)>t\})&\leq& \frac 1t \intrd |f(x)|d\nu(x)\\
		& \leq& \frac 1 t\left(\intrd|f(x)|d\nu_0(x)+\|f\|_{Lip} W_1(\nu,\nu_0)\right)\\
		&\leq& \frac 1t\left(\intrd|f(x)|d\nu_0(x)+\|f\|_{Lip}\sqrt{\theta}\right)\\
		&\leq& 1-\beta, \hfill \text{ for $t$ sufficiently large}.	
		\end{eqnarray*}
		 Thus $\varnb(f)\leq t$, independently of $\nu$. 
	 Similarly $\nu(\{x:\ f(x)<-t\})\leq \beta$ for $t$ sufficiently large, giving us a lower bound $\vmin$. 
		
	Since $\alpha$ ranges over a compact set, and the function $(\alpha,\nu)\mapsto  \alpha+\rrbeta\intrd (f(x)-\alpha)^+\dnu(x)$, being convex in $\alpha$ and linear in $\nu$, is easily seen to be ``convex-concavelike" \cite{borwein_fans_1986}, we may use Ky Fan's minimax theorem. Thus
	\begin{eqnarray*}
		\tr(f) &=& \min_{\vmin\leq \alpha\leq \vmax}\sup_{\nin} \left\{\alpha+\rrbeta\intrd (f(x)-\alpha)^+\dnu(x)\right\}\\
		&=& \mina \sup_{\nin}
		\intrd \alpha+\rrbeta (f(x)-\alpha)^+ d\nu(x) \\
		&=& \mina \left\{\alpha+\rrbeta\sup_{\nin} 
		\intrd   (f(x)-\alpha)^+ d\nu(x)\right\}.
	\end{eqnarray*}

	Applying robust duality Equation (\ref{eq:dual-robust}) to the function $x\mapsto (f(x)-\alpha)^+$, we get 
	\begin{eqnarray*}
		\tr(f) &=& \mina \left\{ \alpha+\rrbeta\infl \left[
		\lambda\theta+\intrd  ((f(x)-\alpha)^+)^{\lambda c} d\nu_0(x)\right] \right\} \\
		&=& \mina\infl \left\{ \alpha+\rrbeta 
		\lambda\theta+\rrbeta\intrd   ((f(x)-\alpha)^+)^{\lambda c} d\nu_0(x) \right\} \\
	\end{eqnarray*}
Writing the minimum as an infimum, and combining the infimums, we get the result.
\end{proof}

We demonstrate Theorem \ref{th:robust-ES} with a simple toy example. 

\begin{example}
Consider a uniform distribution $\nu_0$ on the unit interval $[0,1]$, and let $f(x)=x$ for $x\in\RR$. The minimum in Equation (\ref{eq:dual-robust}) can be obtained in closed form. The first-order conditions for $\alpha$ and $\lambda$ yield the solutions
$\alpha=\beta+\frac{1}{2\lambda}$ and $\lambda=\sqrt{\frac{1-\beta}{2\theta}}$. This corresponds to the transport map 
$$T_{\lambda}(x)=\begin{cases} 
	x+\frac{1}{\lambda} &  \text{ if $\beta<x<1$} \\
	x & \text{ if $0<x\leq \beta$}
\end{cases},$$ and we see that the constraint $\dc(\nu_0,\nu)\leq\theta$ is satisfied, where $\nu=(T_{\lambda})_\# (\nu_0)$, since $\int_0^1 \frac{1}{2}|T_{\lambda}(x)-x|^2 dx=\half\int_{\beta}^1 (\rlambda)^2 dx=\theta$. As could be expected, the robustification procedure shifts the upper tail of the distribution as far as possible while respecting the transport cost constraint. 
\end{example}

The next example calculates the robust, relative to the risk-neutral measure, Expected Shortfall of an unprotected call option, from the point of view of the writer. 

\begin{example}\label{ex:es-call}
Consider a call option on a share that has a risk-neutral distribution given by $\nu_0$, with payoff $f(x)=(x-k)^+$, where $k$ is the strike price. The price of a call option, with strike $k$, will be denoted by $\call(k)$. We can also assume $\alpha\geq 0$, remembering that $\alpha_0=\var(f)$ is a minimizer, since $f\geq 0$, and thus the duality formula will be applied to
$g(x)=(f(x)-\alpha)^+=((x-k)^+-\alpha)^+=(x-(k+\alpha))^+=\max\{x-(k+\alpha),0\}$.

By Theorem \ref{thm:lc} $$g^{\lambda c}(x)=\max\left\{x-(k+\alpha)+\frac{1}{2\lambda},0\right\}=\left(x-(k+\alpha-\frac{1}{2\lambda})\right)^+, \text{ for }\lambda>0.$$ For $\lambda=0$ we clearly have $g^{\lambda c}(x)=\infty$, which can be disregarded in the minimization over $\lambda$.

Since $\nuu$ is a pricing measure, 
$\int 	g^{\lambda c}(x)\dnuu(x)=\call(k+\alpha-\frac{1}{2\lambda})$. By Theorem \ref{th:robust-ES},
\begin{equation}
	\tr(f)=\mina\infl \left\{\alpha+\rrbeta\left[\lambda\theta+\call(k+\alpha-\rtlambda)\right]\right\}.
\end{equation}

The first-order conditions yield
$\lambda=\sqrt{\frac{1-\beta}{2\theta}}$ and $k+\alpha-\rtlambda=\qb$,
where $\qb$ is the $\beta$-quantile of $\nun$, defined by 
$\int_{\qb}^{\infty}\dnuu=1-\beta.$ Therefore
\begin{equation}
	\tr(f)=(\qb-k)+\rrbeta\call(\qb)+\sqrt{\frac{2\theta}{1-\beta}}.
\end{equation}
We compare this with the non-robustified $\cvarnnb(f)$. If $k\leq\qb$ then
\begin{eqnarray*}
	\rho_{\nu_0}(f):=\cvarnnb(f) &=& \rrbeta \int_{\qb}^{\infty} (x-k)d\nu_0(x)\\
	&=& \rrbeta \left(\intqb x-\qb\ d\nun +\intqb \qb-k \ d\nun\right) \\
	&=& \rrbeta\call(\qb)+(\qb-k).
\end{eqnarray*}
If $k>\qb$ then similarly $\rho_{\nu_0}(f)=\rrbeta \call(k)$. 

In contrast to robustifications of ES based on the distribution of $f(X)$,  \cite{wozabal_robustifying_2014}\cite[Table 1]{bartl_computational_2019}, the ``robustification correction" $\tr(f)-\rho_{\nu_0}(f)$ term for Expected Shortfall is a function of $f$, since it depends on $k$. Indeed, 
\begin{equation}
\tr(f)=\rho_{\nu_0}(f)+\sqrt{\frac{2\theta}{1-\beta}}+
\begin{cases}
	\qb-k &  \text{ if } k\leq \qb, \\
	0  & \text{ if } k> \qb.
\end{cases}
\end{equation}
This difference is due to the definition of the risk measure in terms of the payoff $f$ instead of the distribution, and the fact that in the ES maximization process, imagining an adversary that wants to maximize risk, it is only worthwhile to transport mass that is already on the right hand side of the strike price.

In conclusion, we see that robustification adds a term proportional to the square root of the ambiguity radius, and inversely proportional to $1-\beta$. In particular, this quantifies the high distributional sensitivity of an Expected Shortfall calculation for a confidence level $\beta\approx 1$.
\end{example}

\section{Conclusion}
The transform of a payoff $f$ needed for distributional robustification measured with optimal transport of measures under quadratic cost, is tractable if $f$ is convex, especially if it is a maximum of a finite number of affine function. Combining this with the dual representation of Expected Shortfall, it is possible to derive an analytical formula for robust Expected Shortfall, subject to an ambiguity tolerance parameter $\theta$. 
%For a three-asset exposure involving call and put options, distributional robust calculations, subject to a choice of $\theta$, gives a quantitative sense of the distributional risk introduced by these options. 
For the quadratic cost function in optimal transport theory, distributional robust calculations for the Expected Shortfall of convex payoffs are not much more complex than the non-robust calculations.

\section*{Acknowledgement} The work is based on research supported in part by the National Research Foundation of South
Africa (Grant Number 146018). 
%The author also acknowledges the hospitality of the TU Delft and the University of Vienna, and discussions with Antonis Papapantoleon and Daniel Bartl. 

\bibliographystyle{plainurl}
\bibliography{references}

@article{bertsimas_shortfall_2004,
	title = {Shortfall as a risk measure: properties, optimization and applications},
	volume = {28},
	issn = {0165-1889},
	doi = {https://doi.org/10.1016/S0165-1889(03)00109-X},
	number = {7},
	journal = {Journal of Economic Dynamics and Control},
	author = {Bertsimas, Dimitris and Lauprete, Geoffrey J. and Samarov, Alexander},
	year = {2004},
	pages = {1353--1381},
}

@article{borwein_fans_1986,
	title = {On {Fan}'s minimax theorem},
	volume = {34},
	issn = {0025-5610},
	doi = {10.1007/BF01580587},
	number = {2},
	journal = {Mathematical Programming},
	author = {Borwein, J. M. and Zhuang, D.},
	year = {1986},
	mrnumber = {838482},
	pages = {232--234},
}

@article{BreuerCsiszár2016,
	author = {Breuer, Thomas and Csiszár, Imre},
	title = {Measuring Distribution Model Risk},
	journal = {Mathematical Finance},
	volume = {26},
	number = {2},
	pages = {395-411},
	doi = {https://doi.org/10.1111/mafi.12050},
	year = {2016}
}

@article{bartl_computational_2019,
	title = {Computational aspects of robust optimized certainty equivalents and option pricing},
	volume = {2019},
	doi = {10.1111/mafi.12203},
	journal = {Mathematical Finance},
	author = {Bartl, Daniel and Drapeau, Samuel and Tangpi, Ludovic},
	year = {2019},
	pages = {1--23},
}

@book{morini2011understanding,
	title={Understanding and Managing Model Risk: A Practical Guide for Quants, Traders and Validators},
	author={Morini, M.},
	isbn={9780470977743},
	doi={10.1002/9781118467312},
	lccn={2011031397},
	series={The Wiley Finance Series},
	year={2011},
	publisher={Wiley}
}

@misc{feng2018quantifyingmodelriskinherent,
	title={Quantifying the Model Risk Inherent in the Calibration and Recalibration of Option Pricing Models}, 
	author={Yu Feng and Ralph Rudd and Christopher Baker and Qaphela Mashalaba and Melusi Mavuso and Erik Schlögl},
	journal = {arXiv preprint},
	year={2018},
	archivePrefix={arXiv},
	primaryClass={q-fin.RM},
	note={arXiv preprint},
	url={https://arxiv.org/abs/1810.09112}, 
}

@article{delbaen_monetary_2024,
	title = {Monetary utility functions on ${C}_b({X})$ spaces},
	volume = {27},
	doi = {10.1142/S0219024923500334},
	number = {03n04},
	journal = {International Journal of Theoretical and Applied Finance},
	author = {Delbaen, Freddy},
	year = {2024},
	note = {\_eprint: https://doi.org/10.1142/S0219024923500334},
	pages = {2350033},
}

@article{glasserman_robust_2014,
	title = {Robust risk measurement and model risk},
	volume = {14},
	doi = {10.1080/14697688.2013.822989},
	number = {1},
	journal = {Quantitative Finance},
	author = {Glasserman, Paul and Xu, Xingbo},
	year = {2014},
	note = {Publisher: Taylor \& Francis},
	pages = {29--58},
}

@article{rockafellar_optimization_2000,
	title = {Optimization of {C}onditional {V}alue-at-{R}isk},
	volume = {2},
	issn = {14651211},
	doi = {10.21314/JOR.2000.038},
	language = {en},
	number = {3},
	urldate = {2024-09-16},
	journal = {The Journal of Risk},
	author = {Rockafellar, R. Tyrrell and Uryasev, Stanislav},
	year = {2000},
	pages = {21--41},
}

@book{villani_optimal_2009,
	series = {Grundlehren der {Mathematischen} {Wissenschaften}},
	title = {Optimal transport, {Old} and new},
	volume = {338},
	isbn = {978-3-540-71049-3},
	url = {https://doi.org/10.1007/978-3-540-71050-9},
	publisher = {Springer-Verlag, Berlin},
	author = {Villani, Cédric},
	year = {2009},
	doi = {10.1007/978-3-540-71050-9},
	mrnumber = {2459454},
}

@article{wozabal_robustifying_2014,
	title = {Robustifying convex risk measures for linear portfolios: a nonparametric approach},
	volume = {62},
	issn = {0030-364X,1526-5463},
	url = {https://doi.org/10.1287/opre.2014.1323},
	doi = {10.1287/opre.2014.1323},
	number = {6},
	journal = {Operations Research},
	author = {Wozabal, David},
	year = {2014},
	mrnumber = {3294546},
	pages = {1302--1315},
}
\end{document}